\documentclass[12pt,reqno,a4paper]{amsart}
\usepackage{amssymb,amsfonts,amsmath,amsbsy,amsthm,amscd}
\usepackage[
	left = 2.4cm,
	right = 2.4cm,
	top = 2.2 cm,
	bottom = 2.2 cm
]{geometry}
\usepackage[
	pdfencoding=auto,
	colorlinks = true,
	citecolor = blue,
	linkcolor = blue,
	anchorcolor = blue,
	urlcolor = blue
]{hyperref}

\usepackage{amsaddr}
\usepackage{mathtools}

\newcommand{\bb}[1]{\mathbb{#1}}
\newcommand{\bbr}{\bb{R}}

\newcommand{\bv}{\mathbf{v}}

\newcommand{\bn}{\mathbf{n}}

\newcommand{\cE}{\mathcal{E}}

\newcommand{\onehalf}{\frac{1}{2}}
\newcommand{\abs}[1]{\left\vert #1 \right \vert}

\newcommand{\norm}[1]{\left\Vert #1 \right \Vert}
\newcommand{\normm}[1]{\Vert #1 \Vert}
\newcommand{\sqbractau}[1]{\left[ #1 \right]_{t = 0}^{t = \tau}}

\newcommand{\ddt}{\frac{\mathrm{d}}{\mathrm{dt}}}

\newcommand{\intqtau}{\int_{Q_\tau}}
\newcommand{\pt}{\partial_t}

\newcommand{\Div}{\operatorname{div}}
\renewcommand{\d}{\mathrm{d}}
\newcommand{\dx}{\,\mathrm{d} x}
\newcommand{\dt}{\,\mathrm{d} t}
\newcommand{\dxdt}{\,\mathrm{d}x \mathrm{d} t}
\newcommand{\dtau}{\,\mathrm{d} \tau}
\newcommand{\ep}{\varepsilon}

\newcommand{\Ov}[1]{ \overline{ #1 } }

\newcommand{\Grad}{\nabla}
\newcommand{\vr}{\varrho}
\newcommand{\vre}{\varrho_{\ep}}

\newcommand{\vc}[1]{ {\bf #1} }
\newcommand{\tn}[1]{\mathbb{ #1 }}
\newcommand{\R}{\mathbb{ R }}

\numberwithin{equation}{section}

\newcommand{\eps}{\ensuremath{\varepsilon}}

\newtheorem{thm}{Theorem}[section]

\newtheorem{defn}[thm]{Definition}

\newtheorem{rem}[thm]{Remark}

\allowdisplaybreaks[4]

\begin{document}

\title[Low Mach Number Limit for Compressible Two-Phase Flows]{Low Mach Number Limit of a Diffuse Interface Model for Two-Phase Flows of Compressible Viscous Fluids}

\author[H. Abels, Y. Liu, and \v{S}. Ne\v{c}asov\'a]{
	\small
	Helmut Abels$^\ast$, 
	Yadong Liu$^{\dagger\ast}$, and 
	\v{S}\'arka Ne\v{c}asov\'{a}$^\ddagger$
}

\address{
	$^\ast$Fakult\"at f\"ur Mathematik,
	Universit\"at Regensburg,
	93053 Regensburg, Germany
}
\email{helmut.abels@ur.de}

\address{
	$^\dagger$School of Mathematical Sciences, 
	Nanjing Normal University, 
	Nanjing 210023, China
}
\email{ydliu@njnu.edu.cn(yadong.liu@ur.de)}

\address{
	$^\ddagger$Institute of Mathematics,
	Czech Academy of Sciences, 115 67 Praha, Czech Republic
}
\email{matus@math.cas.cz}

\date{\today}

\begin{abstract}
In this paper, we consider a singular limit problem for a diffuse interface model for two immiscible compressible viscous fluids. Via a relative entropy method, we obtain a convergence result for the low Mach number limit to a corresponding system for incompressible fluids in the case of well-prepared initial data and same densities in the limit. 
\end{abstract}



\maketitle

\section{Introduction and Main Result}

Diffuse interface models are an important modeling approach to describe two- or multi-phase flows in fluid mechanics. In comparison with classical sharp interface models they have the theoretical and practical advantage that surfaces separating the fluids do not need to be resolved explicitely. In the case of two fluids the (diffuse) interface is described as the region, where an order parameter, which will be the concentration difference of the two fluids in the following, is not close to one of two values, which describe the presence of only one fluid ($\pm 1$ in the following).

In this contribution we consider the relation between two diffuse interface models for a two-phase flow of viscous Newtonian fluids. The first one is for the case of compressible fluids and leads to the Navier-Stokes/Cahn-Hilliard system for compressible fluids:
\begin{alignat}{2}\label{eq:NSCH1}
  \varrho \partial_t \bv + \varrho \bv\cdot \Grad \bv -\Div \tn{S} + \tfrac1M \Grad p &= -\Div
  \left(\Grad c \otimes \Grad c-\tfrac{|\Grad c|^2}2\tn{I}\right), \\\label{eq:NSCH2}
  \partial_t\varrho +\Div (\varrho \bv) &= 0, \\\label{eq:NSCH3}
  \varrho\partial_t  c + \varrho \bv\cdot \Grad c &= \Delta \mu, \\\label{eq:NSCH4}
  \varrho \mu &= \varrho\tfrac1M \tfrac{\partial f}{\partial c} -\Delta c,
\end{alignat}
in $\Omega\times (0,T)$, where $\Omega\subseteq \R^3$ is a bounded $C^2$-domain and $p = \varrho^2 \frac{\partial f}{\partial \varrho}(\varrho,c)$ and
\begin{align}
  \label{eq:ViscousStress}
  \mathbb{S} &= 2\nu(c)
D\bv + \eta(c) \Div \bv\, \tn{I},\qquad \\\nonumber
D\bv &= \tfrac12 (\Grad \bv + \Grad \bv^T) - \tfrac13 \Div \bv \tn{I}.
\end{align}
Here $c\colon \Omega\times (0,T)\to \R$ describes the concentration difference of two partly miscible compressible fluids, $\rho\colon \Omega\times (0,T)\to [0,\infty)$ is the density of the fluid mixture and $\bv\colon \Omega\times (0,T)\to \R^3$ its (barycentric) mean velocity. Moreover, $\lambda, \eta\colon \R\to [0,\infty)$ are functions describing the shear and bulk viscosity of the mixture, $f\colon [0,\infty)\times \R\to \R$ is a homogeneous free energy density of the mixture and $M>0$ is an analogue of a Mach number. Precise assumptions will be given below. This system is a variant of the model derived by Lowengrub and Truskinovsky~\cite{LowengrubQuasiIncompressible} in a non-dimensionalized form, cf. also~\cite{DiffIntModels}. Here we have set the Reynolds and Peclet number to one for simplifity and the Cahn number proportional to $M$, cf.~\cite[Equations (3.35)]{LowengrubQuasiIncompressible} for the details. We note that in the present variant the total free energy is given by
\begin{equation}\label{eq:FreeEnergy}
  E_{\rm free}(\varrho,c)= \int_\Omega \left(\varrho f(\varrho,c)+\frac12 |\Grad c|^2\right)\dx,
\end{equation}
while in \cite{LowengrubQuasiIncompressible} there is an additional factor $\rho$ in front of $|\nabla c|^2$.
The system is closed by the initial and boundary
conditions
\begin{align}
  \label{m4}
\bv|_{\partial \Omega} &=\Grad c \cdot \vc{n}|_{\partial
    \Omega}= \Grad \mu \cdot \vc{n}|_{\partial \Omega} = 0,\\
  (\bv,c)|_{t=0} &= (\bv_0,c_0),
\end{align}
where $\vc{n}$ is the exterior normal of $\Omega$.
Existence of weak solutions for this system was proved by Feireisl and the first author in \cite{CompressibleNSCH}. This result was extended to the case of certain dynamic boundary conditions by Cherfils et al.~\cite{CherfilsEtAlCompNSCHDynamicBCs}. Existence and uniqueness of strong solutions for this system was shown by Kotschote and Zacher~\cite{KotschoteZacher}, see also~\cite{KotschoteHandbook}. Existence of dissipative martingal solutions of a stochastically perturbed version of this system was shown by Feireisl and Petcu~\cite{FeireislPetcuStochCompNSCH}. In the time-independent, stationary situation existence of weak solutions was shown by Liang and Wang~\cite{LiangWangStatCompNSAC,LiangWangStatCompNSAC2}. An entropy stable finite volume method for the instationary system was proposed by Feireisl, Petcu and She~\cite{FeireislPetcuShe23}, where also existence of weak solutions of the discretized system was shown. Existence of weak solutions for a similar Navier-Stokes/Allen-Cahn system for compressible fluids was shown by 
Feireisl et al.~\cite{FeireislEtAlCompNSAC}. For this system Feireisl, Petcu, and Pra\v{z}\'{a}k~\cite{FereislPetcuPrazak19} studied a relative entropy and obtained results on weak-strong uniqueness and on a low Mach number limit similar to our result in the following.

It is the goal of this contribution to study the low Mach number limit $M\to 0$ for \eqref{eq:NSCH1}-\eqref{eq:NSCH4} and show convergence to solutions of the system
\begin{subequations}
 	\label{eqs:limit-model}
 	\begin{alignat}{2}
 		\label{eqs:limit-v}
 		\partial_t \bv + \bv\cdot \nabla \bv -\Div (2\nu(c)D\bv) + \nabla p & = -\Div (\nabla c\otimes \nabla c),\\
 		\Div \bv & = 0,\\
 		\label{eqs:limit-c}
 		\partial_t c + \bv \cdot \nabla c & = \Delta \mu,\\
 		\label{eqs:limit-mu}
 		\mu & = -\Delta c + G'(c).
 	\end{alignat}
       \end{subequations}
under suitable assumptions and well-prepared initial data. We note that we consider a situation, where the two fluids in limit $M\to 0$ have the same density (or the density difference is neglected). The latter system is known as ``model H'' and is one of the basic diffuse interface models for the two-phase flow of incompressible fluids. It first appeared in Hohenberg and Halperin~\cite{HohenbergHalperin} and was later derived in the framework of rational continuum mechanics by Gurtin et al.~\cite{GurtinTwoPhase}. A first analytic result on existence of strong solutions, if $\Omega=\R^2$ and $G$ is a suitably smooth double well potential was obtained by Starovoitov~\cite{StarovoitovModelH}. More complete results were presented by Boyer~\cite{BoyerModelH} in the case that $\Omega\subset \R^d$ is a periodical channel and a smooth double well potential $G$ and the first author in \cite{ModelH} in the case of a bounded smooth domain and singular double well potential. We refer to Abels, Giorgini, and Garcke~\cite{AGGio} for recent analytic results for an extension of this model to different densities and further references.

We note that, using \eqref{eqs:limit-mu}, one observes that  \eqref{eqs:limit-v} is equivalent to
\begin{equation*}
  \partial_t \bv + \bv\cdot \nabla \bv -\Div (2\nu(c)D\bv) + \nabla p  = \mu \nabla c - \nabla G(c).
\end{equation*}

The mathematical study of the low Mach number limit for systems of equations describing a motion of fluids gets back to the seminal work of Klainerman and Majda \cite{KM1}. Studying various types of singular limits allows us to eliminate unimportant or unwanted modes of the motion as a consequence of scaling and asymptotic analysis. 
The aim of the mathematical analysis of low Mach number limits is to fill up the gap between compressible fluids and their "idealized" incompressible models. There are two ways to introduce the Mach number into the system, which are different from the physical point of view, but  from the mathematical one - completely equivalent. The first approach considers a varying equation of state as well as the transport coefficients see the works of Ebin \cite{EB1},  Schochet \cite{SCH2}. 
The second way is to evaluate qualitatively the incompressibility using the dimensional analysis.  We rewrite our system in the dimensionless form by scaling each variable by its characteristic value, see Klein \cite{Kl}.
The mathematical analysis of singular limits in the frame of strong solutions can be referred to works of Gallagher \cite{Gallag}, Schochet \cite{SCH2}, Danchin \cite {Da} or Hoff \cite{Ho}. 
The seminal works of Lions \cite{LI4} and its extension by Feireisl et al. \cite{FNP} on the existence of global weak solutions in the barotropic case gave a new possibility of a rigorous study of singular limits in the frame of weak solutions, see the works of Desjardins and Grenier \cite{DesGre}, Desjardins, Grenier, Lions and Masmoudi \cite{DGLM}. 

The relative energy inequality was introduced by Dafermos \cite{D} and in the fluid dynamic context was introduced by Germain \cite{G}.
Deriving the relative energy inequality for sufficiently smooth test functions 
and proving  the weak-strong uniqueness it gives us very  powerful and elegant tool for the purpose of measuring the stability
of a solution compared to another solution with better regularity.
This method was developed
by Feireisl, Novotn\'y and co-workers in the framework of singular limits problems (see for example \cite{FJN2012,FN,FN_1,FN2017}
 and references therein).

The structure of this contribution is as follows: In Section~\ref{sec:main} we summarize our assumptions, basic definitions and state our main result on the low Mach number limit. Then in Section~\ref{sec:Proof} we prove the main result with the aid of a relative entropy method.

\subsection*{Notations} In the manuscript, we denote the usual Lebesgue and Sobolev spaces by $ L^p $ and $ W^{k,p} $ respectively for $ 1 \leq p \leq \infty $, $ k \geq 0 $. The corresponding norms are $ \normm{\cdot}_{L^p} $ and $ \normm{\cdot}_{W^{k,p}} $. In particular, we define $ H^k \coloneqq W^{k,2} $. Throughout the paper, the letter $ C $ will indicate a generic positive constant that may change its value from line to line, or even in the same line.

\section{Assumptions and Main Result}\label{sec:main}

We assume that $f$ is given in the form
\begin{equation}
\label{v9}
f(\varrho,c) = f_{\rm e}(\varrho) + M G(c).
\end{equation}
This choice coincides with the assumptions in \cite{CompressibleNSCH} with $ H \equiv 0 $ therein. We only added the factor $M$ in front of $G$, which can be incorporated in $G$.
This yields
\begin{equation}
\label{v10}
p(\varrho,c) = \varrho^2 \frac{\partial f(\varrho, c)}{\partial
\varrho} = p_{\rm e}(\varrho)
,\
f_{\rm e}(\varrho) = \int_1^\varrho \frac{p_e(z)}{z^2} \ {\rm d}z  
\end{equation}
where $p_{\rm e}\in C([0,\infty))\cap C^1(0,\infty)$.
Moreover, it was assumed that
\begin{equation}
\label{v11}
p_{\rm e}(0) = 0,\
\underline{p}_1 \varrho^{\gamma - 1} - \underline{p}_2 \leq
p'_{\rm e}(\varrho) \leq \Ov{p}( 1 + \varrho^{\gamma - 1})
\end{equation}
for a certain $\gamma>\frac32$ and
\begin{equation}
\label{v12}
\begin{gathered}
	G''(c) \geq - \kappa \text{ for some } \kappa \in \bbr, \ 
	\underline{G}_1 c - \underline{G}_2 \leq G'(c) \leq \Ov{G}(1 + c),\\
	\abs{G'(c_1) - G'(c_2)} \leq \overline{G} \abs{c_1 - c_2},\ 
	\abs{G''(c_1) - G''(c_2)} \leq \overline{G} \abs{c_1 - c_2}
\end{gathered}
\end{equation}
for all $c,c_1,c_2\in\R$. Hence \eqref{eq:NSCH1}-\eqref{eq:NSCH4} reduce to
\begin{subequations}
  \label{eqs:scaled-model}
  \begin{alignat}{2}
    \label{eqs:bv-varep}
    \varrho_\eps \partial_t \bv_\eps + \varrho_\eps \bv_\eps \cdot \Grad \bv_\eps -\Div \tn{S}_\eps + \frac1{\eps^2}\Grad (p_e(\varrho_\eps) - p_e(1)) & = \vr_\varepsilon \mu_\varepsilon \nabla c_\varepsilon - \vr_\varepsilon G'(c_\varepsilon) \nabla c_\varepsilon, \\
    \label{eqs:vr-varep}
    \partial_t\varrho_\eps +\Div (\varrho_\eps \bv_\eps) & = 0, \\
    \label{eqs:c-varep}
    \varrho_\eps\partial_t  c_\eps + \varrho_\eps \bv_\eps\cdot \Grad c_\eps & = \Delta \mu_\eps, \\
    \label{eqs:mu-varep}
    \varrho_\eps \mu_\eps & = \vr_\eps G'(c_\eps) -\Delta c_\eps,
	\end{alignat}
	subject to the boundary conditions
	\begin{equation}
		\bv_\varepsilon|_{\partial \Omega} =\nabla c_\varepsilon \cdot \bn|_{\partial
		\Omega}= \nabla \mu_\varepsilon \cdot \bn|_{\partial \Omega} = 0.
            \end{equation}
 \end{subequations}

Let us recall the definition of weak solutions in the sense of \cite[Theorem~1.2]{CompressibleNSCH} (with $H\equiv 0$ there):
\begin{defn}
  Let $T>0$, $Q_T=\Omega\times (0,T)$, $\vr_{0,\eps}\in L^\gamma(\Omega)$ with $\vr_{0,\eps}\geq 0$ almost everywhere, and $\vc{m}_{0,\eps}\colon \Omega\to \R^3$ be measurable such that $\vr_{0,\eps}^{-1}|\vc{m}_{0,\eps}|^2\in L^1(\Omega)$. Then $\varrho_\eps\in
  L^\infty(0,T;L^\gamma(\Omega))$ with $\varrho_\eps\geq 0$, $\bv_\eps\in L^2(0,T;H^1(\Omega; \R^3))$, $c_\eps\in L^\infty(0,T;H^1(\Omega))$ are a weak solution of \eqref{eqs:scaled-model} if the following holds true:
  \begin{enumerate}
  \item For every $\boldsymbol \varphi \in \mathcal{D}(\Omega \times (0,T); \R^3)$
    \begin{eqnarray}\nonumber
      \lefteqn{-\int_{Q_T}\Big( \vr_\eps \bv_\eps \cdot\partial_t \boldsymbol{\varphi} + \left(\vr_\eps \bv_\eps
        \otimes \bv_\eps + \tfrac1M p_e(\vr_\eps)\, \tn{I} -  \tn{S}_\eps\right):
        \Grad \boldsymbol\varphi\Big)\dxdt 
      }\\\label{a5}
      &=& \int_{Q_T} \Big( (\Grad c_\eps \otimes \Grad c_\eps): \Grad \boldsymbol
        \varphi - \frac{|\Grad c_\eps|^2}{2} \Div \boldsymbol \varphi  \Big)\dxdt,\qquad 
    \end{eqnarray}
    where $\mathbb{S}_\eps = 2\nu(c_\eps)D\bv_\eps + \eta(c_\eps) \Div \bv_\eps\, \tn{I}$.
  \item $\varrho_\eps$ is a renormalized solution of \eqref{eqs:vr-varep} in the sense of
    DiPerna and Lions~\cite{DiPernaLions}, i.e.,
    \begin{equation}\label{a2}
    \int_{Q_T} \Big( \vr_\eps B(\vr_\eps) \partial_t \varphi + \vr_\eps B(\vr_\eps) \bv_\eps
      \cdot \Grad \varphi - b(\vr_\eps) \Div \bv_\eps \ \varphi \Big) \dxdt = 0
    \end{equation}
    for any test function $\varphi \in \mathcal{D}( \Ov{\Omega}\times (0,T))$,
    and any
    \begin{equation}\label{a3}
    B(\vr_\eps) = B(1) + \int_1^{\vr_\eps} \frac{b(z)}{z^2} \,\mathrm d z,
    \end{equation}
    where $b\in C^0([0,\infty))$ is a bounded function.
  \item For every $\varphi \in \mathcal{D}(\Omega\times (0,T) )$
    \begin{equation}\label{m3}
      \int_{Q_T}\left(\varrho_\eps c_\eps\,\partial_t \varphi+ \varrho_\eps c_\eps \bv_\eps\cdot \Grad
        \varphi\right)\dxdt =\int_{Q_T} \Grad \mu_\eps\cdot\Grad \varphi\dxdt
    \end{equation}
    and
    \begin{equation}\label{cp}
      \int_{Q_T}\varrho_\eps \mu_\eps\varphi\dxdt = \int_{Q_T}\left(\varrho_\eps G'(c_\eps)\varphi +
        \Grad c_\eps\cdot \Grad \varphi\right)\dxdt.
    \end{equation}
  \item The
\emph{energy inequality}
\begin{equation}\label{a6}
  E(t)+ \int_{Q_\tau} \left( \tn{S}_\eps: \Grad
  \bv_\eps + | \Grad \mu_\eps |^2 \right) \, \d x \d \tau
  \leq E(s)
\end{equation}
holds for almost every $0\leq s\leq T$ including $s=0$ and all $t\in [s,T]$, where
\begin{align}\label{eq:Energy1}
  E(t) &= \int_\Omega  \varrho_\eps(t) \frac{|\bv_\eps(t)|^2}2 \dx + E_{\rm free}(\varrho_\eps(t),c_\eps(t)),\\
  \label{eq:Energy2}
  E(0)= E_0 &= \int_\Omega \varrho^{-1}_{0,\eps} \frac{|\vc{m}_{0,\eps}|^2}2 \dx + E_{\rm free}(\varrho_{0,\eps},c_{0,\eps}).
\end{align}
\item $\varrho_\eps,\varrho_\eps\bv_\eps,c_\eps$ are weakly continuous with respect to $t\in[0,T]$ with values in $L^1(\Omega)$ and
  $\varrho_\eps|_{t=0}=\varrho_{0,\eps}$, $\varrho_\eps\bv_\eps|_{t=0}=\vc{m}_{0,\eps}$, $c_\eps|_{t=0}=c_{0,\eps}$.
\end{enumerate}
\end{defn}

We note that existence of solutions follows from \cite[Theorem~1.2]{CompressibleNSCH}.

Our main result in this contribution is:
\begin{thm}
	\label{thm:main}
	Let {$\gamma\geq \frac{12}5$}, $\Omega\subseteq\R^3$ be a bounded domain with $C^2$-boundary, $T>0$, $M=\eps^2$, $\eps>0$, and let $(\bv,c,\mu)$ be a (sufficiently) smooth solution of \eqref{eqs:limit-model}. 
	Moreover, we assume that $\vr_{0,\eps}\in L^\gamma(\Omega)$, $\bv_{0,\eps}\in L^2(\Omega)^3$, $c_{0,\eps}\in H^1(\Omega)$ are given such that $\vr_{0, \varepsilon} \coloneqq 1 + \varepsilon \vr_{0, \varepsilon}^{(1)}$ and
	\begin{equation}
		\label{eqs:intial-well}
		\normm{\vr_{0, \varepsilon}^{(1)}}_{L^\infty(\Omega)}
		+ \norm{\bv_{0, \varepsilon}}_{L^2(\Omega)}
		+ \norm{\nabla c_{0, \varepsilon}}_{L^2(\Omega)}
		\leq C.
	\end{equation}
	and $\vr_{0,\eps}^{(1)}\to_{\eps\to 0} 0$ in $L^\infty(\Omega)$, $\bv_{0,\eps}\to_{\eps\to 0} \bv_0$ in $L^2(\Omega)^3$, $c_{0,\eps}\to_{\eps\to 0} c_0$ in $H^1(\Omega)$. Then
	\begin{equation*}
		\vr_\eps(t,\cdot)\to_{\eps\to 0}1 \text{ in }L^1(\Omega),\quad \bv_\eps(t,\cdot)\to_{\eps\to 0}\bv\text{ in } L^2(\Omega)^3,\quad c_\eps(t,\cdot)\to_{\eps\to 0} c(t,\cdot) \text{ in } H^1(\Omega)
	\end{equation*}
	uniformly in $t\in [0,T]$.
\end{thm}

\begin{rem}
		Here the restriction of $ \gamma \geq \frac{12}{5} $ comes essentially from addressing the nonconvex part of the potential $ G $, cf. \eqref{eqs:mix-relative-2}. If one considers a convex potential $ G $, we may relax it to $ \gamma \geq 2 $, which is due to the convection $ \vre \mu_\varepsilon \bv \cdot \nabla c_\varepsilon $, cf. \eqref{eqs:gamma-2}.
\end{rem}

\section{Low Mach Limit}\label{sec:Proof}
\subsection{Relative energy inequality}
      To compare \eqref{eqs:scaled-model} and \eqref{eqs:limit-model}, we proceed with the so-called \textit{relative energy (entropy)}:
\begin{align*}
	\cE\big(\vr_\varepsilon, \bv_\varepsilon, c_\varepsilon | 1, \bv, c\big)
	& \coloneqq \int_\Omega \left[
		\onehalf \vr_\varepsilon \abs{\bv_\varepsilon - \bv}^2
		+ \frac{1}{\varepsilon^2} \big(F_e(\vr_\varepsilon) - F_e'(1)(\vr_\varepsilon - 1) - F_e(1)\big)
	\right]	\dx \\
	& \quad + \int_\Omega \left[
		\onehalf \abs{\nabla c_\varepsilon - \nabla c}^2 
		+ \vr_\varepsilon (G(c_\eps) - G'(c)(c_\varepsilon - c) - G(c))
	\right] \dx,
\end{align*}
where $	F_e(\vr_\varepsilon)= \vr_\varepsilon f_e(\vr_\varepsilon)$.

By the weak formulation of \eqref{eqs:bv-varep} for $ \bv_\eps $, i.e., \eqref{a5}, with $ \bv $ as the test function we obtain for every $\tau \in [0,T]$
\begin{equation}
	\label{eqs:rho-v_e-v}
	\begin{aligned}
		- \sqbractau{\int_\Omega \vr_\eps \bv_\eps \bv \dx}
		& = - \intqtau \vr_\eps \bv_\eps \cdot \pt \bv \dxdt \\
		& \quad 
		- \intqtau \big(
		\vr_\eps \bv_\eps \otimes \bv_\eps
		- 2 \nu(c_\varepsilon) D\bv_\eps
		\big) : \nabla \bv \dxdt \\
		& \quad 
		- \intqtau \vr_\varepsilon \mu_\varepsilon \nabla c_\varepsilon \cdot \bv \dxdt
		+ \intqtau \vr_\varepsilon \bv \cdot \nabla c_\varepsilon G'(c_\varepsilon) \dxdt.
	\end{aligned}
      \end{equation}
Let $ \tfrac{1}{2} \abs{\bv}^2 $ be the test function in the weak formulation of continuity equation \eqref{eqs:vr-varep}. This yields
\begin{equation}
	\label{eqs:rho-v^2}
	\sqbractau{\int_\Omega \frac{\vr_\eps}{2} \abs{\bv}^2 \dx}
	= \intqtau (\vr_\eps \bv \cdot \pt \bv + \vr_\eps \bv_\eps \cdot \nabla \bv \cdot \bv) \dxdt
      \end{equation}
      for every $\tau\in [0,T]$.
Summing \eqref{eqs:rho-v_e-v}, \eqref{eqs:rho-v^2} and the energy inequality of weak solutions, one ends up with
\begin{align}
	\nonumber
	& \sqbractau{\int_\Omega \frac{\vr_\eps}{2} \abs{\bv_\eps - \bv}^2 \dx} \\
	\nonumber
	& \quad + \sqbractau{\int_\Omega \Big(\vr_\varepsilon G(c_\eps) + \onehalf \abs{\nabla c_\eps}^2\Big) \dx}
	+ \sqbractau{\int_\Omega \frac{1}{\varepsilon^2} F_e (\vr_\eps) \dx} \\
	\label{eqs:rho-v_e-v^2}
	& \quad + \intqtau 2 \nu(c_\varepsilon) \abs{D \bv_\varepsilon}^2 \dxdt
	+ \intqtau \abs{\nabla \mu_\varepsilon}^2 \dxdt \\
	\nonumber
	& \leq - \intqtau \vr_\eps (\bv_\eps - \bv) \cdot (\pt \bv + \bv \cdot \nabla \bv)\dxdt
	- \intqtau \vr_\eps (\bv_\eps - \bv) \cdot \nabla \bv \cdot (\bv - \bv_\varepsilon) \dxdt \\
	\nonumber
	& \quad + \intqtau 2 \nu(c_\varepsilon) D\bv_\eps : \nabla \bv \dxdt
	- \intqtau \vr_\varepsilon \mu_\varepsilon \nabla c_\varepsilon \cdot \bv \dxdt
	+ \intqtau \vr_\varepsilon \bv \cdot \nabla c_\varepsilon G'(c_\varepsilon) \dxdt.
\end{align}
In view of the weak formulation of \eqref{eqs:c-varep}, we have
\begin{equation}
	\label{eqs:rho-c-mu}
	\sqbractau{\int_\Omega \vr_\varepsilon c_\varepsilon \mu \dx}
	- \intqtau (\vr_\varepsilon c_\varepsilon \pt \mu + \vr_\varepsilon c_\varepsilon \bv_\varepsilon \cdot \nabla \mu) \dxdt
	= - \intqtau \nabla \mu_\varepsilon \cdot \nabla \mu \dxdt.
\end{equation}
Moreover, with $ \mu = G'(c) - \Delta c $ one gets
\begin{equation}
	\label{eqs:pt-rho-c-mu}
	\intqtau \vr_\varepsilon c_\varepsilon \pt \mu \dxdt
	= \intqtau \vr_\varepsilon c_\varepsilon \pt (G'(c) - \Delta c) \dxdt.
\end{equation}
Adding \eqref{eqs:rho-c-mu} and \eqref{eqs:pt-rho-c-mu} and integration by parts gives
\begin{align}
	\nonumber
	& \sqbractau{- \int_\Omega \vr_\varepsilon c_\varepsilon G'(c) - \nabla c_\varepsilon \cdot \nabla c \dx} \\
	\nonumber
	& = \sqbractau{\int_\Omega (\vr_\varepsilon - 1) c_\varepsilon \Delta c \dx}
	- \intqtau \vr_\varepsilon c_\varepsilon \bv_\varepsilon \cdot \nabla \mu \dxdt \\
	\label{eqs:-rho-ce-G'}
	& \quad + \intqtau \nabla \mu_\varepsilon \cdot \nabla \mu \dxdt
	- \intqtau \vr_\varepsilon c_\varepsilon (G''(c) \pt c - \Delta \pt c) \dxdt.
\end{align}
Direct calculations yield
\begin{align*}
	\sqbractau{\int_\Omega \vr_\varepsilon c G'(c) \dx}
	& = \sqbractau{\int_\Omega (\vr_\varepsilon - 1) c G'(c) \dx}
	+ \int_0^\tau \ddt \int_\Omega cG'(c) \dxdt \\
	& = \sqbractau{\int_\Omega (\vr_\varepsilon - 1) c G'(c) \dx}
	+ \intqtau \pt c G'(c) \dxdt
	+ \intqtau c G''(c) \pt c \dxdt.
\end{align*}
Then we have
\begin{equation}
	\label{eqs:rho-c-G'}
	\begin{aligned}
		\sqbractau{\int_\Omega \vr_\varepsilon c G'(c) - \vr_\varepsilon G(c) \dx}
		& = \sqbractau{\int_\Omega (\vr_\varepsilon - 1) c G'(c) \dx} \\
		& \quad + \sqbractau{\int_\Omega (1 - \vr_\varepsilon) G(c) \dx}
		+ \intqtau c G''(c) \pt c \dxdt.
	\end{aligned}
\end{equation}
It follows from the strong formulation of \eqref{eqs:limit-model} that
\begin{equation}
	\label{eqs:c^2}
	\begin{aligned}
		\sqbractau{\int_\Omega \frac{1}{2}\abs{\nabla c}^2 \dx}
		& = 
		- \intqtau G'(c) \pt c \dxdt 
		- \intqtau \abs{\nabla \mu}^2 \dxdt 
		- \intqtau \bv \cdot \nabla c \mu \dxdt.
	\end{aligned}
\end{equation}
Now we summarize from \eqref{eqs:rho-v_e-v^2}, \eqref{eqs:-rho-ce-G'}, \eqref{eqs:rho-c-G'}, \eqref{eqs:c^2} that
\begin{align*}
	& \sqbractau{\int_\Omega \frac{\vr_\eps}{2} \abs{\bv_\eps - \bv}^2 \dx} 
	+ \sqbractau{\int_\Omega \frac{1}{\varepsilon^2} F_e (\vr_\eps) \dx} \\
	& \quad + \sqbractau{\int_\Omega \Big(\vr_\varepsilon (G(c_\eps) - G'(c)(c_\varepsilon - c) - G(c)) + \onehalf \abs{\nabla c_\eps - \nabla c}^2\Big) \dx} \\
	& \quad + \intqtau 2 \nu(c_\varepsilon) D \bv_\varepsilon : (D \bv_\varepsilon - D \bv) \dxdt
	+ \intqtau \abs{\nabla \mu_\varepsilon - \nabla \mu}^2 \dxdt \\
	& \leq - \intqtau \vr_\eps (\bv_\eps - \bv) \cdot (\pt \bv + \bv \cdot \nabla \bv)\dxdt
	- \intqtau \vr_\eps (\bv_\eps - \bv) \cdot \nabla \bv \cdot (\bv - \bv_\varepsilon) \dxdt \\
	& \quad + \sqbractau{\int_\Omega (\vr_\varepsilon - 1) c_\varepsilon \Delta c \dx}
	+ \sqbractau{\int_\Omega (\vr_\varepsilon - 1) c G'(c) \dx} 
	+ \sqbractau{\int_\Omega (1 - \vr_\varepsilon) G(c) \dx} \\
	& \quad 
	- \intqtau \vr_\varepsilon \mu_\varepsilon \nabla c_\varepsilon \cdot \bv \dxdt
	+ \intqtau \vr_\varepsilon \bv \cdot \nabla c_\varepsilon G'(c_\varepsilon) \dxdt \\
	& \quad
	- \intqtau \vr_\varepsilon c_\varepsilon \bv_\varepsilon \cdot \nabla \mu \dxdt 
	- \intqtau \vr_\varepsilon c_\varepsilon (G''(c) \pt c - \Delta \pt c) \dxdt \\
	& \quad 
	+ \intqtau c G''(c) \pt c \dxdt
	- \intqtau G'(c) \pt c \dxdt 
	- \intqtau \bv \cdot \nabla c \mu \dxdt.
\end{align*}
In view of the strong formulation of \eqref{eqs:limit-v},
\begin{align*}
	& - \intqtau \vr_\eps (\bv_\eps - \bv) \cdot (\pt \bv + \bv \cdot \nabla \bv)\dxdt \\
	& = - \intqtau \vr_\eps (\bv_\eps - \bv) \cdot \big(\Div(2 \nu(c) D\bv) - \nabla \pi + \mu \nabla c - \nabla G(c)\big)\dxdt\\
	\nonumber
	& = - \intqtau (\vr_\eps - 1) (\bv_\eps - \bv) \cdot \big(\Div(2 \nu(c) D\bv) + \mu \nabla c\big)\dxdt \\
	\nonumber
	& \quad 
	- \intqtau \vr_\varepsilon (\bv_\varepsilon - \bv) \cdot \nabla (\pi + G(c)) \dxdt \\
	\nonumber
	& \quad 
	- \intqtau (\bv_\eps - \bv) \cdot \big(\Div(2 \nu(c) D\bv)  \big)\dxdt - \intqtau (\bv_\eps - \bv) \cdot \mu \nabla c \dxdt.
\end{align*}
By integration by parts, one obtains
\begin{align}
	\nonumber
	& - \intqtau (\bv_\eps - \bv) \cdot \big(\Div(2 \nu(c) D\bv) \big)\dxdt \\
	\nonumber
	& = \intqtau (D\bv_\varepsilon - D\bv) : (2 \nu(c) D\bv) \dxdt \\
	\nonumber
	& = \intqtau (D\bv_\varepsilon - D\bv) : (2 \nu(c_\varepsilon) D\bv) \dxdt
	+ \intqtau (D\bv_\varepsilon - D\bv) : (2 (\nu(c) - \nu(c_\varepsilon)) D\bv) \dxdt
\end{align}
Then adding all together with \eqref{eqs:vr-varep} multiplied by $ F_e'(\vr_\varepsilon) $ entails that
\begin{align}
	\nonumber
	& \Big[\cE\big(\vr_\varepsilon, \bv_\varepsilon, c_\varepsilon | 1, \bv, c\big)\Big]_{t = 0}^{t = \tau} \\
	\nonumber
	& \quad + \intqtau 2 \nu(c_\varepsilon) \abs{D \bv_\varepsilon - D \bv}^2 \dxdt
	+ \intqtau \abs{\nabla \mu_\varepsilon - \nabla \mu}^2 \dxdt \\
	\nonumber
	& \leq \sqbractau{\int_\Omega (\vr_\varepsilon - 1) c_\varepsilon \Delta c \dx}
	+ \sqbractau{\int_\Omega (\vr_\varepsilon - 1) c G'(c) \dx}
	+ \sqbractau{\int_\Omega (1 - \vr_\varepsilon) G(c) \dx} \\
	\nonumber
	& \quad 
	- \intqtau (\vr_\eps - 1) (\bv_\eps - \bv) \cdot \Div(2 \nu(c) D\bv)\dxdt \\
	\nonumber
	& \quad - \intqtau \vr_\eps (\bv_\eps - \bv) \cdot \nabla \bv \cdot (\bv - \bv_\varepsilon) \dxdt
	- \intqtau \vr_\varepsilon (\bv_\varepsilon - \bv) \cdot \nabla (\pi + G(c)) \dxdt \\
	\nonumber
	& \quad - \intqtau \vr_\varepsilon (\bv_\eps - \bv) \cdot \nabla(c \mu) \dxdt 
	- \intqtau \vr_\varepsilon (\bv_\eps - \bv) \cdot \nabla \mu (c_\varepsilon - c) \dxdt \\
	\nonumber
	& \quad - \intqtau (\vr_\varepsilon - 1) \bv \cdot \nabla \mu c_\varepsilon \dxdt
	- \intqtau \bv \cdot \nabla(c_\varepsilon - c) (\mu_\varepsilon - \mu) \dxdt \\
	\nonumber
	& \quad
	- \intqtau (\vr_\varepsilon - 1) \mu_\varepsilon \nabla c_\varepsilon \cdot \bv \dxdt
	+ \intqtau (\vr_\varepsilon - 1) \bv \cdot \nabla c_\varepsilon G'(c_\varepsilon) \dxdt \\
	\nonumber
	& \quad
	+ \intqtau (D\bv_\varepsilon - D\bv) : (2 (\nu(c) - \nu(c_\varepsilon)) D\bv) \dxdt \\
	\nonumber
	& \quad 
	+ \intqtau \vr_\varepsilon \pt c \big(G'(c_\varepsilon) - (c_\varepsilon - c) G''(c) - G'(c)\big) \dxdt \\
	\nonumber
	& \quad 
	- \intqtau (\vr_\varepsilon - 1) \big(G''(c) c \pt c - G'(c) \pt c\big) \dxdt \\
	\label{eqs:relative-energy-inequlaity}
	& \quad - \intqtau (\vr_\varepsilon - 1) c_\varepsilon \Delta \pt c \dxdt
	- \intqtau (\vr_\varepsilon - 1) \mu_\varepsilon \pt c \dxdt,
\end{align}
where we used
\begin{align*}
	\intqtau \vr_\varepsilon \bv \cdot \nabla c_\varepsilon G'(c_\varepsilon) \dxdt
	& = \intqtau (\vr_\varepsilon - 1) \bv \cdot \nabla c_\varepsilon G'(c_\varepsilon) \dxdt
	+ \intqtau \bv \cdot \nabla G(c_\varepsilon) \dxdt \\
	& = \intqtau (\vr_\varepsilon - 1) \bv \cdot \nabla c_\varepsilon G'(c_\varepsilon) \dxdt.
\end{align*}

Concerning the potential part of $G$, here due to the assumption \eqref{v12} we employ the decomposition of $ G $ such that $ G(c) = G_0(c) + G_1(c) $ with $ G_1(c) = - \kappa \frac{c^2}{2} $ for $ \kappa > 0 $ where $ G_0(c) $ is convex. Then we have
\begin{equation*}
		\int_\Omega \vr_\varepsilon (G_0(c_\eps) - G_0'(c)(c_\varepsilon - c) - G_0(c)) \dx \geq 0,
\end{equation*}
and
\begin{align}
	\label{eqs:mix-relative-3}
	- \sqbractau{\int_\Omega \vr_\varepsilon (G_1(c_\eps) - G_1'(c)(c_\varepsilon - c) - G_1(c)) \dx} 
	= \sqbractau{\int_\Omega \vr_\varepsilon  \frac{\kappa}{2} (c_\varepsilon - c)^2 \dx}.
\end{align}
As $ G_1(c) $ is a nonconvex part, in the following we would like to justify the following identity to ensure a suitable relative energy inequality:
	\begin{align}
		\nonumber
		& \sqbractau{\int_\Omega   \frac{\vr_\varepsilon}{2} (c_\varepsilon - c)^2 \dx} \\
		\nonumber
		& \quad
		= - \intqtau (\nabla \mu_\varepsilon - \nabla \mu) \cdot (\nabla c_\varepsilon - \nabla c) \dxdt
		- \intqtau (\vre - 1) (c_\varepsilon - c) \pt c \dxdt
		\\
		\label{eqs:mix-relative-2}
		& \qquad
		- \intqtau (c_\varepsilon - c) (\bv_\varepsilon - \bv) \cdot \nabla c \dxdt
		- \intqtau (\vre - 1)(c_\varepsilon - c) \bv_\varepsilon \cdot \nabla c \dxdt.
	\end{align}
        
We give an essential claim for the justification.

\noindent
\textbf{Claim:} It holds that
\begin{align}
	\label{eqs:rho-c-ve-2}
	& \sqbractau{\int_\Omega \vr_\varepsilon \frac{c_\eps^2}2 \dx}  =
	-\intqtau \nabla \mu_\varepsilon \cdot \nabla c_\varepsilon  \dxdt.
\end{align}
\emph{Proof of the claim.} Let $ 0 \leq  t \leq t + h \leq T $. Integrating \eqref{eqs:vr-varep} over $[t,t+h]$ in its weak formulation (using a standard approximation argument) and testing with $\frac12c_\eps (t+h)c_\eps(t)$ yields
\begin{alignat*}{1}
  \int_\Omega \frac{(\vr_\eps(t+h)-\vr_\eps(t))c_\eps (t+h)c_\eps(t)}{2h} \dx  = \int_\Omega\frac1{2h} \int_t^{t+h}  \vr_\eps(\tau) \bv_\eps(\tau)\dtau \cdot \nabla (c_\eps(t+h)c_\eps(t)) \dx
\end{alignat*}
Similarly, integrating \eqref{eqs:c-varep} on $[t,t+h]$ in its weak formulation (using a standard approximation argument) and testing with $\frac12(c_\eps (t+h)+c_\eps(t))$ yields
\begin{align*}
  &\int_\Omega \frac{(\vr_\eps(t+h)c_\eps(t+h)-\vr_\eps(t)c_\eps(t))(c_\eps (t+h)+c_\eps(t))}{2h} \dx \\
  &= \int_\Omega\frac1{2h} \int_t^{t+h}  (\vr_\eps(\tau)c_\eps(\tau) \bv_\eps(\tau)-\nabla \mu_\eps(\tau)) \dtau \cdot \nabla (c_\eps(t+h)+c_\eps(t)) \dx
\end{align*}
Now subtracting the first from the second identity gives
\begin{alignat*}{1}
  &\int_\Omega \frac{\vr_\eps(t+h)c_\eps^2(t+h)-\vr_\eps(t)c_\eps^2(t)}{2h} \dx\\
  &= \int_\Omega\frac1{2h} \int_t^{t+h}  (\vr_\eps(\tau)c_\eps(\tau) \bv_\eps(\tau)-\nabla \mu_\eps(\tau)) \dtau \cdot \nabla (c_\eps(t+h)+c_\eps(t)) \dx\\
  &\quad -\int_\Omega\frac1{2h} \int_t^{t+h}  \vr_\eps(\tau) \bv_\eps(\tau) \dtau \cdot \nabla (c_\eps(t+h)c_\eps(t)) \dx.
\end{alignat*}
To pass to the limit $h\to 0+$ in the first term on the right-hand side it is essential that $\nabla (c_\eps(.+h)+c_\eps)\in L^2(0,T;L^p(\Omega))$ and $\frac1{2h} \int_t^{t+h}  \vr_\eps(\tau)c_\eps(\tau) \bv_\eps(\tau)) \dtau\in L^2(0,T;L^{p'}(\Omega))$ are bounded, where $\frac1p =\frac1\gamma -\frac16$, cf.\ \eqref{eqs:nablac-L2Lp} below. This is the case if $\frac1\gamma + \frac16+\frac16 \leq 1-\frac1\gamma +\frac16$, which is equivalent to $\gamma\geq \frac{12}5$. This is where we need an extra restriction on $ \gamma $. The same estimates can be applied for the second term on the right-hand side. Hence we can pass to the limit $h\to 0 +$ and obtain
\begin{align*}
  \ddt \int_\Omega \vr_\eps \frac{c_\eps^2}2 \dx = -\int_\Omega \nabla \mu_\eps \cdot \nabla c_\eps \dx  \quad \text{ in }\mathcal{D}'(0,T),
\end{align*}
which yields the claim by the fundamental theorem for Sobolev functions and integrating over $ (0,T) $. \hfil\qed

Taking $ c $ as the test function in the weak formulation of \eqref{eqs:c-varep} and employing the strong formulation of \eqref{eqs:limit-c}, we find
\begin{align}
	\nonumber
	\sqbractau{\int_\Omega \vr_\varepsilon c_\eps c \dx}
	& = \intqtau (- \nabla \mu_\varepsilon + \vre c_\varepsilon \bv_\varepsilon) \cdot \nabla c \dxdt
	+ \int_\Omega \vre c_\varepsilon \pt c \dx \\
	\label{eqs:rho-c-ve-c}
	& = \intqtau (\vre - 1) c_\varepsilon \bv_\varepsilon \cdot \nabla c \dxdt
	+ \intqtau c_\varepsilon (\bv_\varepsilon - \bv) \cdot \nabla c \dxdt \\
	\nonumber
	& \quad
	- \intqtau \nabla \mu_\varepsilon \cdot \nabla c \dxdt
	- \intqtau \nabla \mu \cdot \nabla c_\varepsilon \dxdt 
	+ \int_\Omega (\vre - 1) c_\varepsilon \pt c \dxdt.
\end{align}

Similarly, taking $ \frac{c^2}{2} $ as the test function in the weak formulation of the continuity equation \eqref{eqs:vr-varep}, together with \eqref{eqs:limit-c}, yields
\begin{align}
	\nonumber
	\sqbractau{\int_\Omega \vr_\varepsilon \frac{c^2}{2} \dx}
	& = \intqtau \vre c \bv_\varepsilon \cdot \nabla c \dxdt
	+ \intqtau \vre c \pt c \dxdt \\
	\nonumber
	& = \intqtau (\vre - 1) c \bv_\varepsilon \cdot \nabla c \dxdt
	+ \intqtau c (\bv_\varepsilon - \bv) \cdot \nabla c \dxdt \\
	\label{eqs:rho-c-2}
	& \quad
	- \intqtau \nabla \mu \cdot \nabla c \dxdt
	+ \int_\Omega (\vre - 1) c \pt c \dxdt.
\end{align}
Summing \eqref{eqs:rho-c-ve-2}, \eqref{eqs:rho-c-2}, and subtracting \eqref{eqs:rho-c-ve-c} from the resulting equation entail the desired identity \eqref{eqs:mix-relative-2}.

Now we define a modified relative energy $ \widetilde{\cE} $ by eliminating the nonconvex part of the chemical potential
\begin{align*}
	\widetilde{\cE}\big(\vr_\varepsilon, \bv_\varepsilon, c_\varepsilon | 1, \bv, c\big)
	\coloneqq \cE\big(\vr_\varepsilon, \bv_\varepsilon, c_\varepsilon | 1, \bv, c\big)
	- \int_\Omega \vr_\varepsilon (G_1(c_\eps) - G_1'(c)(c_\varepsilon - c) - G_1(c)) \dx.
\end{align*}
Adding \eqref{eqs:relative-energy-inequlaity} and \eqref{eqs:mix-relative-2}, one obtains \textit{relative energy inequality}
\begin{align}
	\nonumber
	& \Big[\widetilde{\cE}\big(\vr_\varepsilon, \bv_\varepsilon, c_\varepsilon | 1, \bv, c\big)\Big]_{t = 0}^{t = \tau}
  \\
	\nonumber
	& \quad + \intqtau 2 \nu(c_\varepsilon) \abs{D \bv_\varepsilon - D \bv}^2 \dxdt
	+ \intqtau \abs{\nabla \mu_\varepsilon - \nabla \mu}^2 \dxdt \\
	\nonumber
	& \leq \sqbractau{\int_\Omega (\vr_\varepsilon - 1) c_\varepsilon \Delta c \dx}
	+ \sqbractau{\int_\Omega (\vr_\varepsilon - 1) c G'(c) \dx}
	+ \sqbractau{\int_\Omega (1 - \vr_\varepsilon) G(c) \dx} \\
	\nonumber
	& \quad 
	- \kappa \intqtau (\nabla \mu_\varepsilon - \nabla \mu) \cdot (\nabla c_\varepsilon - \nabla c) \dxdt
		- \kappa \intqtau (\vre - 1) (c_\varepsilon - c) \pt c \dxdt
		 \\
	\nonumber
	& \quad - \kappa \intqtau (c_\varepsilon - c) (\bv_\varepsilon - \bv) \cdot \nabla c \dxdt
		- \kappa \intqtau (\vre - 1)(c_\varepsilon - c) \bv_\varepsilon \cdot \nabla c \dxdt \\
	\nonumber
	& \quad 
	- \intqtau (\vr_\eps - 1) (\bv_\eps - \bv) \cdot \Div(2 \nu(c) D\bv)\dxdt \\
	\nonumber
	& \quad - \intqtau \vr_\eps (\bv_\eps - \bv) \cdot \nabla \bv \cdot (\bv - \bv_\varepsilon) \dxdt
	- \intqtau \vr_\varepsilon (\bv_\varepsilon - \bv) \cdot \nabla (\pi + G(c)) \dxdt \\
	\nonumber
	& \quad - \intqtau \vr_\varepsilon (\bv_\eps - \bv) \cdot \nabla(c \mu) \dxdt 
	- \intqtau \vr_\varepsilon (\bv_\eps - \bv) \cdot \nabla \mu (c_\varepsilon - c) \dxdt \\
	\nonumber
	& \quad - \intqtau (\vr_\varepsilon - 1) \bv \cdot \nabla \mu c_\varepsilon \dxdt
	- \intqtau \bv \cdot \nabla(c_\varepsilon - c) (\mu_\varepsilon - \mu) \dxdt \\
	\nonumber
	& \quad
	- \intqtau (\vr_\varepsilon - 1) \mu_\varepsilon \nabla c_\varepsilon \cdot \bv \dxdt
	+ \intqtau (\vr_\varepsilon - 1) \bv \cdot \nabla c_\varepsilon G'(c_\varepsilon) \dxdt \\
	\nonumber
	& \quad
	+ \intqtau (D\bv_\varepsilon - D\bv) : (2 (\nu(c) - \nu(c_\varepsilon)) D\bv) \dxdt \\
	\nonumber
	& \quad 
	+ \intqtau \vr_\varepsilon \pt c \big(G'(c_\varepsilon) - (c_\varepsilon - c) G''(c) - G'(c)\big) \dxdt \\
	\nonumber
	& \quad 
	- \intqtau (\vr_\varepsilon - 1) \big(G''(c) c \pt c - G'(c) \pt c\big) \dxdt \\
	\label{eqs:m-relative-energy-inequlaity}
	& \quad - \intqtau (\vr_\varepsilon - 1) c_\varepsilon \Delta \pt c \dxdt
	- \intqtau (\vr_\varepsilon - 1) \mu_\varepsilon \pt c \dxdt.
\end{align}

\subsection{Uniform estimates}
Let $ \bv = 0 $ and $ c = 1 $ in \eqref{eqs:m-relative-energy-inequlaity}. Then one obtains
\begin{equation*}
	\Big[\widetilde\cE\big(\vr_\varepsilon, \bv_\varepsilon, c_\varepsilon | 1, 0, 1\big)\Big]_{t = 0}^{t = \tau} + \intqtau 2 \nu(c_\varepsilon) \abs{D \bv_\varepsilon}^2 \dxdt 
	+ \intqtau \abs{\nabla \mu_\varepsilon}^2 \dxdt \leq C
      \end{equation*}
      for every $\tau \in [0,T]$.
In a similar way as in \cite{FJN2012,FereislPetcuPrazak19}, we obtain the uniform estimates
\begin{gather}
	\operatorname*{esssup}_{t \in (0,T)} \norm{\sqrt{\vr_\varepsilon} \bv_\varepsilon}_{L^2} \leq C, \\
	\label{eqs:uniform-vr^2}
	\operatorname*{esssup}_{t \in (0,T)} \int_{\Omega \cap \{1/2 \leq \vr_\varepsilon \leq 2\}} \abs{\frac{\vr_\varepsilon - 1}{\varepsilon}}^2 \dx \leq C, \\
	\label{eqs:uniform-vr^gamma}
	 \operatorname*{esssup}_{t \in (0,T)} \int_{\Omega \setminus \{1/2 \leq \vr_\varepsilon \leq 2\}} (1 + \abs{\vr_\varepsilon}^\gamma) \dx \leq \varepsilon^2 C, \\
	 \label{eqs:uniform-nabla-c}
	 \operatorname*{esssup}_{t \in (0,T)} \norm{\nabla c_\varepsilon}_{L^2(\Omega)}
	 \leq C, \\
	 \label{eqs:uniform-nabla-mu}
	 \int_0^T \norm{\nabla \mu_\varepsilon}_{L^2(\Omega)}^2 \dt \leq C, 
\end{gather}
where $ C > 0 $ depends on the bounds for the initial data. Moreover, via Korn's inequality (cf. \cite[Theorem 11.21]{FN2017}) and $ \nu_* \leq \nu(c_\varepsilon) \leq \nu^* $, one has
\begin{equation*}
	\int_0^T \norm{\nabla \bv_\varepsilon}_{L^2(\Omega)}^2 \dt \leq C.
\end{equation*}
In view of a generalized Korn--Poincar\'e inequality (cf. \cite[Theorem 11.23]{FN2017}), we obtain
\begin{equation}
	\label{eqs:uniform-v}
	\int_0^T \norm{\bv_\varepsilon}_{W^{1,2}(\Omega)}^2 \dt \leq C.
\end{equation}
Incorporating with \eqref{eqs:uniform-nabla-c} and the conservation of $ \vr_\varepsilon c_\varepsilon $, i.e., for a.e. $\tau \in (0,T)$,
\begin{equation*}
	\int_\Omega \vr_\varepsilon c_\varepsilon(\tau) \dx
	= \int_\Omega \vr_0 c_0 \dx,
\end{equation*}
proceeding in a similar way as in \cite[Lemma 2.1]{CompressibleNSCH} yields
\begin{equation}
	\label{eqs:uniform-c}
	\operatorname*{esssup}_{t \in (0,T)} \norm{c_\varepsilon(t)}_{W^{1,2}(\Omega)}
	\leq C.
\end{equation}

Moreover, it follows from \eqref{v12}, \eqref{eqs:mu-varep}, \eqref{eqs:uniform-c} that
\begin{equation*}
	\left|\int_\Omega \vre \mu_\varepsilon \dx\right|
	= \left|\int_\Omega \vre G'(c_\varepsilon) \dx\right|
	\leq C \norm{\vre}_{L^\gamma(\Omega)} \bigg(1 + \norm{c_\varepsilon}_{L^{\frac{\gamma}{\gamma - 1}}(\Omega)}\bigg) \leq C,
\end{equation*}
which, in accordance with \eqref{eqs:uniform-nabla-mu}, implies
\begin{equation}
	\label{eqs:uniform-mu}
	\int_0^T \norm{\mu_\varepsilon(t)}_{W^{1,2}(\Omega)}^2\,dt
	\leq C.
\end{equation}
With Sobolev embedding in 3D, we have $ \mu_\varepsilon \in L^2(0,T; L^6(\Omega)) $. Combining with the fact $ \vre \in L^\infty(0,T; L^\gamma(\Omega)) $, one obtains $ \vre \mu_\varepsilon \in L^2(0,T; L^q(\Omega)) $ uniformly, with $ \frac{1}{q} = \frac{1}{\gamma} + \frac{1}{6} $. Then by means of the elliptic estimates of $ c_\eps $ in \eqref{eqs:mu-varep}, namely,
\begin{equation*}
	- \Delta c_\eps
	= \varrho_\eps \mu_\eps - \vr_\eps G'(c_\eps),
\end{equation*}
we get
\begin{equation}
	\label{eqs:uniform-c-elliptic}
	c_\eps \in L^2(0,T; W^{2,q}(\Omega))
\end{equation}
for all $ 1 < q < 6 $ satisfying $ \frac{1}{q} = \frac{1}{\gamma} + \frac{1}{6} $.

\subsection{Incompressible limit}
Now we are in the position to control the right-hand side terms of \eqref{eqs:relative-energy-inequlaity} and derive the desired limit passage. First, 
\begin{equation*}
	\intqtau \vr_\eps (\bv_\eps - \bv) \cdot \nabla \bv \cdot (\bv - \bv_\varepsilon) \dxdt
	\leq \int_0^\tau \norm{\nabla \bv}_{L^\infty(\Omega)} \norm{\sqrt{\vr_\varepsilon} (\bv_\varepsilon - \bv)}_{L^2(\Omega)}^2 \dt \leq C \int_0^\tau \widetilde{\cE}(t) \dt.
\end{equation*}
By $ \nu_* \leq \nu \leq \nu^* $, the Lipschitz continuous of $ \nu(c) $, and Young's inequality, we have
\begin{align*}
	& \intqtau (D\bv_\varepsilon - D\bv) : (2 (\nu(c) - \nu(c_\varepsilon)) D\bv) \dxdt \\
	& \leq \onehalf \intqtau \nu(c_\varepsilon) \abs{D\bv_\varepsilon - D\bv}^2 \dxdt
	+ C(\nu_*^{-1}) \int_0^\tau \norm{D\bv}_{L^\infty(\Omega)}^2 \norm{c_\varepsilon - c}_{L^2(\Omega)}^2 \dt \\
	& \leq \frac{1}{2} \intqtau \nu(c_\varepsilon) \abs{D \bv_\varepsilon - D \bv}^2 \dxdt
	+ C(\nu_*^{-1}) \int_0^\tau \widetilde{\cE}(t) \dt.
\end{align*}
Moreover,
\begin{align*}
	& \intqtau \vr_\varepsilon (\bv_\eps - \bv) \cdot \nabla \mu (c_\varepsilon - c) \dxdt \\
	& \leq \int_0^\tau \norm{\vr_\varepsilon}_{L^\gamma(\Omega)} \norm{\nabla \mu}_{L^\infty(\Omega)} \norm{\bv_\varepsilon - \bv}_{L^6(\Omega)} \norm{c_\varepsilon - c}_{L^6(\Omega)} \dt \\
	& \leq \frac{1}{2} \intqtau \nu(c_\varepsilon) \abs{D \bv_\varepsilon - D \bv}^2 \dxdt
	+ C(\nu_*^{-1}) \int_0^\tau \widetilde{\cE}(t) \dt,
\end{align*}
for all $ \gamma > 3/2 $, where we used the energy bounded of $ \vr_\varepsilon $, the Sobolev embedding of $ W^{1,2} \hookrightarrow L^6 $ in three dimensions and the Poincar\'e inequality. Analogously, it follows
\begin{align*}
	- \intqtau \bv \cdot \nabla(c_\varepsilon - c) (\mu_\varepsilon - \mu) \dxdt
	& = \intqtau \bv \cdot \nabla(\mu_\varepsilon - \mu) (c_\varepsilon - c) \dxdt \\
	& \leq \int_0^\tau \norm{\bv}_{L^\infty(\Omega)}^2 \norm{c_\varepsilon - c}_{L^2(\Omega)}^2 \dt 
	+ \frac{1}{4} \intqtau \abs{\nabla \mu_\varepsilon - \nabla \mu}^2 \dxdt \\
	& \leq \frac{1}{4} \intqtau \abs{\nabla \mu_\varepsilon - \nabla \mu}^2 \dxdt
	+ C \int_0^\tau \widetilde{\cE}(t) \dt,
\end{align*}
and
\begin{align*}
	& - \kappa \intqtau (\nabla \mu_\varepsilon - \nabla \mu) \cdot (\nabla c_\varepsilon - \nabla c) \dxdt
	- \kappa \intqtau (c_\varepsilon - c) (\bv_\varepsilon - \bv) \cdot \nabla c \dxdt \\
	& \quad \leq \frac{1}{4} \intqtau \abs{\nabla \mu_\varepsilon - \nabla \mu}^2 \dxdt 
	+ \frac{1}{2} \intqtau \nu(c_\varepsilon) \abs{D \bv_\varepsilon - D \bv}^2 \dxdt
	+ C(\nu_*^{-1}) \int_0^\tau \widetilde{\cE}(t) \dt.
\end{align*}
By direct calculations and weak formulation of continuity equation for $ \vr_\varepsilon $,
\begin{align}
	\nonumber
	& \intqtau \vr_\varepsilon (\bv_\varepsilon - \bv) \cdot \nabla (p + c \mu + G(c)) \dxdt \\
	\nonumber
	& = - \varepsilon \intqtau \frac{\vr_\varepsilon - 1}{\varepsilon} \pt (p + c \mu + G(c)) \dxdt \\
	\label{eqs:rho-v-v}
	& \quad + \varepsilon \sqbractau{\int_\Omega \frac{\vr_\varepsilon - 1}{\varepsilon} (p + c \mu + G(c)) \dx}
	- \varepsilon \intqtau \frac{\vr_\varepsilon - 1}{\varepsilon} \bv \cdot \nabla(p + c \mu + G(c)) \dxdt.
\end{align}
For sufficiently smooth $ (\bv,p,c,\mu) $, it follows from \eqref{eqs:uniform-vr^2} and the H\"older inequality that \eqref{eqs:rho-v-v} is controlled by
\begin{equation*}
	\intqtau \vr_\varepsilon (\bv_\varepsilon - \bv) \cdot \nabla (p + c \mu+G(c)) \dxdt
	\leq \varepsilon C,
\end{equation*}
where $ C $ depends on the initial data and $ (\bv,p,c,\mu) $, but is independent of $ \varepsilon > 0 $. 

Concerning the terms associated with $ (\vr_\varepsilon - 1) $, it follows from \eqref{eqs:uniform-vr^2} and
\eqref{eqs:uniform-vr^gamma} that
\begin{equation*}
	\intqtau (\vr_\varepsilon - 1) f \dxdt
	\leq \varepsilon C,
\end{equation*}
for all $ f \in L^1(0,T; L^2(\Omega) \cap L^{\frac{\gamma}{\gamma - 1}}(\Omega)) $ with $ \gamma > \frac{3}{2} $. 
Similarly,
\begin{align*}
	& \sqbractau{\int_\Omega (\vr_\varepsilon - 1) c_\varepsilon \Delta c \dx}
	+ \sqbractau{\int_\Omega (\vr_\varepsilon - 1) c G'(c) \dx}\\
	& \quad 
	+ \sqbractau{\int_\Omega (\vr_\varepsilon - 1) c G'(c) \dx}
	+ \sqbractau{\int_\Omega (1 - \vr_\varepsilon) G(c) \dx}
	\leq \varepsilon C.
\end{align*}
However, for the term 
$ \intqtau (\vr_\varepsilon - 1) \mu_\varepsilon \nabla c_\varepsilon \cdot \bv \dxdt $,
we know from \eqref{eqs:uniform-c} and \eqref{eqs:uniform-mu} that $ \mu_\varepsilon \nabla c_\varepsilon \in L^2(0,T; L^{\frac{3}{2}}(\Omega)) $, which is not sufficient for all $ \gamma > \frac{3}{2} $.  By \eqref{eqs:uniform-c-elliptic} and Sobolev embedding $ W^{2,q} \hookrightarrow W^{1,p} $ with $ \frac{1}{p} = \frac{1}{q} - \frac{1}{3} = \frac{1}{\gamma} - \frac{1}{6} $, one obtains
\begin{equation}
	\label{eqs:nablac-L2Lp}
	\nabla c_\varepsilon \in L^2(0,T; L^p(\Omega))
\end{equation}
for all $ p $ satisfying $ \frac{1}{p} = \frac{1}{\gamma} - \frac{1}{6} $. Then 
\begin{equation}
	\label{eqs:gamma-2}
	\intqtau (\vr_\varepsilon - 1) \mu_\varepsilon \nabla c_\varepsilon \cdot \bv \dxdt
	\leq \int_0^\tau\normm{\vr_\varepsilon - 1}_{L^\gamma(\Omega)} \normm{\mu_\varepsilon}_{L^6(\Omega)} \normm{\nabla c_\varepsilon}_{L^p(\Omega)} \normm{\bv}_{L^\infty(\Omega)}\,dt
	\leq \varepsilon C
\end{equation}
for $ 1 \geq \frac{1}{\gamma} + \frac{1}{6} + \frac{1}{p} = \frac{2}{\gamma} $, which holds for all $ \gamma \geq 2 $. 

Furthermore, 
\begin{align*}
	& \intqtau \vr_\varepsilon \pt c \big(G'(c_\varepsilon) - (c_\varepsilon - c) G''(c) - G'(c)\big) \dxdt \\
	& \leq C \intqtau \abs{c_\varepsilon - c}^2 \dxdt 
	+ \intqtau (\vr_\varepsilon - 1) \pt c \big(G'(c_\varepsilon) - (c_\varepsilon - c) G''(c) - G'(c)\big) \dxdt \\
	& \leq C \int_0^\tau \widetilde{\cE}(t) \dt + \varepsilon C. 
\end{align*}

Collecting all the estimates above, we then obtain a Gronwall's type inequality:
\begin{equation*}
	\Big[\widetilde{\cE}\big(\vr_\varepsilon, \bv_\varepsilon, c_\varepsilon | 1, \bv, c\big)\Big]_{t = 0}^{t = \tau}
	\leq C \int_0^\tau \widetilde{\cE}\big(\vr_\varepsilon, \bv_\varepsilon, c_\varepsilon | 1, \bv, c\big)(t) \dt + \varepsilon C
      \end{equation*}
      for all $\tau\in (0,T)$,
which yields
\begin{equation*}
	\widetilde{\cE}\big(\vr_\varepsilon, \bv_\varepsilon, c_\varepsilon | 1, \bv, c\big)(\tau)
	\leq C\big(\widetilde{\cE}\big(\vr_{0, \varepsilon}, \bv_{0, \varepsilon}, c_{0, \varepsilon} | 1, \bv_0, c_0\big) + \varepsilon\big) e^{\tau}
      \end{equation*}
      for all $\tau\in (0,T)$.
If additionally one has for the initial data 
\begin{alignat*}{3}
	\bv_{0, \varepsilon} & \to \bv_0, && \text{ in } L^2(\Omega), \\
	\vr_{0, \varepsilon}^{(1)} & \to 0, && \text{ in } L^\infty(\Omega), \\
	\nabla c_{0, \varepsilon} & \to \nabla c_0, && \text{ in } L^2(\Omega),
\end{alignat*}
as $\varepsilon \to 0$,
one concludes the low Mach number limit immediately, which finishes the proof of Theorem \ref{thm:main}.

\section*{Acknowledgement} 

This work was initiated when H. Abels and Y. Liu visited the Institue of Mathematics of the Czech Academy of Science in Prague during November 2022. The hospitality is gratefully appreciated.

The work of \v S. Ne\v casov\' a  has been supported by the Czech Science Foundation (GA\v CR) through projects project 22-01591S and also by  Praemium Academiae of \v S. Ne\v casov\' a.

\section*{Compliance with Ethical Standards}
\subsection*{Date avability}
Data sharing not applicable to this article as no datasets were generated during the current study.
\subsection*{Conflict of interest}
The authors declare that there are no conflicts of interest.


\begin{thebibliography}{10}

\bibitem{ModelH}
{\sc H.~Abels}, {\em On a diffuse interface model for two-phase flows of
  viscous, incompressible fluids with matched densities}, Arch. Rat. Mech.
  Anal., 194 (2009), pp.~463--506.

\bibitem{CompressibleNSCH}
{\sc H.~Abels and E.~Feireisl}, {\em On a diffuse interface model for a
  two-phase flow of compressible viscous fluids}, Indiana Univ. Math. J., 57
  (2008), pp.~659--698.

\bibitem{AGGio}
{\sc H.~Abels, H.~Garcke, and A.~Giorgini}, {\em Global regularity and
  asymptotic stabilization for the incompressible
  {N}avier–{S}tokes-{C}ahn–{H}illiard model with unmatched densities},
  Math. Ann., 55pp.,  (2023).

\bibitem{DiffIntModels}
{\sc D.~M. Anderson, G.~B. McFadden, and A.~A. Wheeler}, {\em Diffuse-interface
  methods in fluid mechanics}, in Annual review of fluid mechanics, Vol. 30,
  vol.~30 of Annu. Rev. Fluid Mech., Annual Reviews, Palo Alto, CA, 1998,
  pp.~139--165.

\bibitem{BoyerModelH}
{\sc F.~Boyer}, {\em Mathematical study of multi-phase flow under shear through
  order parameter formulation}, Asymptot. Anal., 20 (1999), pp.~175--212.

\bibitem{CherfilsEtAlCompNSCHDynamicBCs}
{\sc L.~Cherfils, E.~Feireisl, M.~Mich\'{a}lek, A.~Miranville, M.~Petcu, and
  D.~Pra\v{z}\'{a}k}, {\em The compressible {N}avier-{S}tokes-{C}ahn-{H}illiard
  equations with dynamic boundary conditions}, Math. Models Methods Appl. Sci.,
  29 (2019), pp.~2557--2584.

\bibitem{D}
{\sc C.~M. Dafermos}, {\em The second law of thermodynamics and stability},
  Arch. Rational Mech. Anal., 70 (1979), pp.~167--179.

\bibitem{Da}
{\sc R.~Danchin}, {\em Zero mach number limit for compressible flows with
  periodic boundary conditions}, Amer. J. Math., 124 (2002), pp.~1153--1219.

\bibitem{DesGre}
{\sc B.~Desjardins and E.~Grenier}, {\em Low {M}ach number limit of viscous
  compressible flows in the whole space}, R. Soc. Lond. Proc. Ser. A Math.
  Phys. Eng. Sci., 455 (1999), pp.~2271--2279.

\bibitem{DGLM}
{\sc B.~Desjardins, E.~Grenier, P.-L. Lions, and N.~Masmoudi}, {\em
  Incompressible limit for solutions of the isentropic {N}avier-{S}tokes
  equations with {D}irichlet boundary conditions}, J. Math. Pures Appl. (9), 78
  (1999), pp.~461--471.

\bibitem{DiPernaLions}
{\sc R.~DiPerna and P.~Lions}, {\em {Ordinary differential equations, transport
  theory and Sobolev spaces.}}, Invent. Math., 98 (1989), pp.~511--547.

\bibitem{EB1}
{\sc D.~B. Ebin}, {\em The motion of slightly compressible fluids viewed as a
  motion with strong constraining force}, Ann. Math., 105 (1977), pp.~141--200.

\bibitem{FJN2012}
{\sc E.~Feireisl, B.~J. Jin, and A.~Novotn\'{y}}, {\em Relative entropies,
  suitable weak solutions, and weak-strong uniqueness for the compressible
  {N}avier-{S}tokes system}, J. Math. Fluid Mech., 14 (2012), pp.~717--730.

\bibitem{FN}
{\sc E.~Feireisl and A.~Novotn\'{y}}, {\em Singular limits in thermodynamics of
  viscous fluids}, Advances in Mathematical Fluid Mechanics, Birkh\"{a}user
  Verlag, Basel, 2009.

\bibitem{FN_1}
\leavevmode\vrule height 2pt depth -1.6pt width 23pt, {\em Weak-strong
  uniqueness property for the full {N}avier-{S}tokes-{F}ourier system}, Arch.
  Ration. Mech. Anal., 204 (2012), pp.~683--706.

\bibitem{FN2017}
\leavevmode\vrule height 2pt depth -1.6pt width 23pt, {\em Singular limits in
  thermodynamics of viscous fluids}, Advances in Mathematical Fluid Mechanics,
  Birkh\"{a}user/Springer, Cham, second~ed., 2017.

\bibitem{FNP}
{\sc E.~Feireisl, A.~Novotn\'{y}, and H.~Petzeltov\'{a}}, {\em On the existence
  of globally defined weak solutions to the {N}avier-{S}tokes equations}, J.
  Math. Fluid Mech., 3 (2001), pp.~358--392.

\bibitem{FeireislPetcuStochCompNSCH}
{\sc E.~Feireisl and M.~Petcu}, {\em A diffuse interface model of a two-phase
  flow with thermal fluctuations}, Appl. Math. Optim., 83 (2021), pp.~531--563.

\bibitem{FereislPetcuPrazak19}
{\sc E.~Feireisl, M.~Petcu, and D.~Pra\v{z}\'{a}k}, {\em Relative energy
  approach to a diffuse interface model of a compressible two-phase flow},
  Math. Methods Appl. Sci., 42 (2019), pp.~1465--1479.

\bibitem{FeireislPetcuShe23}
{\sc E.~Feireisl, M.~Petcu, and B.~She}, {\em An entropy stable finite volume
  method for a compressible two phase model}, Appl. Math., 68 (2023),
  pp.~467--483.

\bibitem{FeireislEtAlCompNSAC}
{\sc E.~Feireisl, H.~Petzeltov\'{a}, E.~Rocca, and G.~Schimperna}, {\em
  Analysis of a phase-field model for two-phase compressible fluids}, Math.
  Models Methods Appl. Sci., 20 (2010), pp.~1129--1160.

\bibitem{Gallag}
{\sc I.~Gallagher}, {\em R{\'e}sultats r{\'e}cents sur la limite
  incompressible}, Ast{\'e}risque,S{\'e}minaire Bourbaki, 299 (2005),
  pp.~29--57.

\bibitem{G}
{\sc P.~Germain}, {\em Weak-strong uniqueness for the isentropic compressible
  {N}avier-{S}tokes system}, J. Math. Fluid Mech., 13 (2011), pp.~137--146.

\bibitem{GurtinTwoPhase}
{\sc M.~E. Gurtin, D.~Polignone, and J.~Vi{\~n}als}, {\em Two-phase binary
  fluids and immiscible fluids described by an order parameter}, Math. Models
  Methods Appl. Sci., 6 (1996), pp.~815--831.

\bibitem{Ho}
{\sc D.~Hoff}, {\em Dynamics of singularity surfaces for compressible viscous
  flows in two space dimensions}, Comm. Pure Appl. Math., 55 (2002),
  pp.~1365--1407.

\bibitem{HohenbergHalperin}
{\sc P.~Hohenberg and B.~Halperin}, {\em Theory of dynamic critical
  phenomena.}, Rev. Mod. Phys., 49 (1977), pp.~435--479.

\bibitem{KM1}
{\sc S.~Klainerman and A.~Majda}, {\em Singular limits of quasilinear
  hyperbolic systems with large parameters and the incompressible limit of
  compressible fluids}, Comm. Pure Appl. Math., 34 (1981), pp.~481--524.

\bibitem{Kl}
{\sc R.~Klein, N.~Botta, T.~Schneider, C.~D. Munz, S.~Roller, A.~Meister,
  L.~Hoffmann, and T.~Sonar}, {\em Asymptotic adaptive methods for multi-scale
  problems in fluid mechanics}, vol.~39, 2001, pp.~261--343.
\newblock Special issue on practical asymptotics.

\bibitem{KotschoteHandbook}
{\sc M.~Kotschote}, {\em Local and global existence of strong solutions for the
  compressible {N}avier-{S}tokes equations near equilibria via the maximal
  regularity}, in Handbook of mathematical analysis in mechanics of viscous
  fluids, Springer, Cham, 2018, pp.~1905--1946.

\bibitem{KotschoteZacher}
{\sc M.~Kotschote and R.~Zacher}, {\em Strong solutions in the dynamical theory
  of compressible fluid mixtures}, Math. Models Methods Appl. Sci., 25 (2015),
  pp.~1217--1256.

\bibitem{LiangWangStatCompNSAC}
{\sc Z.~Liang and D.~Wang}, {\em Stationary {C}ahn-{H}illiard-{N}avier-{S}tokes
  equations for the diffuse interface model of compressible flows}, Math.
  Models Methods Appl. Sci., 30 (2020), pp.~2445--2486.

\bibitem{LiangWangStatCompNSAC2}
\leavevmode\vrule height 2pt depth -1.6pt width 23pt, {\em Weak solutions to
  the stationary {C}ahn-{H}illiard/{N}avier-{S}tokes equations for compressible
  fluids}, J. Nonlinear Sci., 32 (2022), pp.~Paper No. 41, 25.

\bibitem{LI4}
{\sc P.-L. Lions}, {\em Mathematical topics in fluid mechanics. {V}ol. 2},
  vol.~10 of Oxford Lecture Series in Mathematics and its Applications, The
  Clarendon Press, Oxford University Press, New York, 1998.
\newblock Compressible models, Oxford Science Publications.

\bibitem{LowengrubQuasiIncompressible}
{\sc J.~Lowengrub and L.~Truskinovsky}, {\em Quasi-incompressible
  {C}ahn-{H}illiard fluids and topological transitions}, R. Soc. Lond. Proc.
  Ser. A Math. Phys. Eng. Sci., 454 (1998), pp.~2617--2654.

\bibitem{SCH2}
{\sc S.~Schochet}, {\em The mathematical theory of low {M}ach number flows},
  M2AN Math. Model Numer. anal., 39 (2005), pp.~441--458.

\bibitem{StarovoitovModelH}
{\sc V.~N. Starovo{\u\i}tov}, {\em On the motion of a two-component fluid in
  the presence of capillary forces}, Mat. Zametki, 62 (1997), pp.~293--305.

\end{thebibliography}

\def\cprime{$'$} \def\ocirc#1{\ifmmode\setbox0=\hbox{$#1$}\dimen0=\ht0
  \advance\dimen0 by1pt\rlap{\hbox to\wd0{\hss\raise\dimen0
  \hbox{\hskip.2em$\scriptscriptstyle\circ$}\hss}}#1\else {\accent"17 #1}\fi}

\end{document}